\documentclass[12pt]{article}
\usepackage{mathrsfs}
\usepackage{amsfonts}
\usepackage{amssymb}
\usepackage{amsthm}
\usepackage{amsmath}
\usepackage{latexsym}
\usepackage{color}

\allowdisplaybreaks[4]
\newtheorem{coro}{\qquad Corollary}
\newtheorem{pron}{\qquad Proposition}

\begin{document}

\title{\large\bf On the almost decrease of a subexponential density
\thanks{Research supported by the National Science Foundation of China
(No. 11071182 {\&} No. 71171177), the project of the
key research base of human and social
science (Statistics, Finance) for colleges in Zhejiang Province (Grant No. of
Academic Education of Zhejiang,~2008-255).}}
\author { Tao Jiang$^{1)}$ \qquad Yuebao Wang$^{2)}$ \thanks{Corresponding author. \qquad Telephone: +86 512 67422726. \qquad Fax: +86 512 65112637. \quad E-mail: ybwang@suda.edu.cn (Y. Wang)}
\qquad Zhaolei Cui$^{3)}$
\qquad\\{\small\it 1). School of Statistics and Mathematics, Zhejiang Gongshang University, P. R. China, 310018 }\\
{\small\it 2). School of Mathematical Sciences, Soochow University, P. R. China, 215006 } \\
{\small\it 3). School of Mathematics and Statistics, Changshu Institute of Technology, P. R. China, 215500 }
}
\date{}

\maketitle

\begin{abstract}
For a subexponential density, so far, there has been no positive conclusion or counter example to show whether it is almost decreasing. In this paper, a subexponential density supported on $\mathbb{R}^+\cup\{0\}$ without the almost decrease is constructed by a little skillful method. The density is a positive piecewise linear function with a more normal shape. Correspondingly, there exists a local subexponential distribution which is not locally almost decreasing. Based on an example of Cline \cite{C1986}, some similar results are also obtained for the long-tailed density excluding the subexponential density and the local long-tailed distribution excluding the local subexponential distribution. Finally, the paper shows that, for the local subexponentiality of a distribution supported on $\mathbb{R}$, the local almost decreasing condition is necessary in some sense.
\medskip

{\it Keywords:} almost decrease; subexponential density; local subexponential distribution\\

{\it 2000 Mathematics Subject Classification:} Primary 60G50, 60F99, 60E05
\end{abstract}

%%%%%%%%%%%%%%%%%%%%%%%%%%%%%%%%%%%%%%%%%Section 1%%%%%%%%%%%%%%%%%%%%%%
\section{Introduction}

In this paper, let $F$ be an absolutely continuous distribution supported on $\mathbb{R}^+\cup\{0\}$ with density $f$, and all limit relationships be for $x\to\infty$, unless otherwise stated. In addition, we write $a(x)\sim b(x)$ for two positive functions $a(\cdot)$ and $b(\cdot)$ whenever $\lim \frac{a(x)}{b(x)}=1$.

We say that a density $f$ belongs to the long-tailed density class, denoted by $f\in\mathcal{L}_0$, if $f(x)>0$ eventually, and for each constant $t\in\mathbb{R}$,
$$f(x+t)\sim f(x).$$
We say that a density $f$ belongs to the subexponential density class, denoted by $f\in\mathcal{S}_0$, if $f\in\mathcal{L}_{0}$ and
$$f^{\otimes2}(x):=\int_0^xf(x-y)f(y)dy\sim 2f(x).$$

The study of subexponential density is an important part of the heavy-tailed distribution theory. It has important applications in many fields of application probability, such as risk theory, queuing system, branching process, and so on. The study can go back to Chover et al. \cite{CNW1973}. Up to now, the more and more related studies and applications have come out, see, for example, Sgibnev \cite{S1996}, Kl\"{u}ppelberg \cite{K1989}, Asmussen et al. \cite{AFK2003}, Wang and Wang \cite{WW2006}, Korshunov \cite{K2006}, Foss et al. \cite{FKZ2013}, and Watanabe and Yamamuro \cite{WY2017}. However, there are still some interesting problems which are worthy of research.

We know that, for a density $f$, though $\int_0^\infty f(y)dy=1$, it does not have to go to zero. If the density $f\in\mathcal{L}_0$, then $f(x)\to0$ in some way, for example, $f$ is almost decreasing. In the terminology of Bingham et al. \cite{BGT1987}, a density $f$ is called almost decreasing, that is almost not increasing, if there exists a constant $x_0\in\mathbb{R}^+\cup\{0\}$ such that $f(x)>0$ for all $x\ge x_0$ and
$$\sup_{x_0\le x\le y<\infty}\frac{f(y)}{f(x)}=:C=C(f,x_0)<\infty.$$
In particular, if $C=1$, then the density is monotonically decreasing. Clearly, a lot of densities are almost decreasing. Of course, there are some densities which are not almost decreasing. Thus, a problem is put forward naturally.

\textbf{Problem 1.1.} Is each density in the class $\mathcal{S}_{0}$ almost decreasing  ?

The problem is closely related to the local asymptotics of a distribution and other research objects, see Problem 1.2, Corollary 1.1, Problem 3.1 and Proposition 3.1 below. So far, we have not found any positive conclusions or counterexamples on this basic problem in distribution theory. Now, we give a negative answer to the problem.
\begin{pron}\label{pron101}
There exists a density $f\in\mathcal{S}_{0}$ without the almost decrease.
\end{pron}

Based on an example of density in Cline \cite{C1986}, we have a corresponding conclusion for another density class.
\begin{pron}\label{pron102}
In the class $\mathcal{L}_{0}\setminus\mathcal{S}_{0}$, there is a density which is not almost decreasing.
\end{pron}

%Because there is a close relationship between density in the class $\mathcal{L}_{0}$ or $\mathcal{S}_{0}$ and corresponding local distribution,
Correspondingly, we can also discuss the locally almost decrease of a distribution. To this end, we first introduce some concepts and notations about the local distribution class.

%In this subsection, we give some concepts of local heavy tailed distributions.
For some constant $d\in\mathbb{R}^+\cup\{\infty\}$ and distribution $F$, we denote $F(x+\Delta_d):=F(x,x+d]$ when $d\in\mathbb{R}^+$ and $\overline{F}(x):=F(x,\infty)$ when $d=\infty$, where $x\in\mathbb{R}$.

We say that a distribution $F$ belongs to the local long-tailed distribution class $\mathcal{L}_{\Delta_d}$, if for some $d\in\mathbb{R}^+\cup\{\infty\},F(x+\Delta_d)>0$ eventually, and for each constant $t\in\mathbb{R}^+$ it holds uniformly for all $\mid s\mid\le t$ that
$$F(x+s+\Delta_d)\sim F(x+\Delta_d).$$
Clearly, the distribution $F$ is heavy-tailed, that is $\int_{0}^{\infty}e^{\varepsilon y}F(dy)=\infty$ holds for each constant $\varepsilon\in\mathbb{R}^+$.
Further, if a distribution $F$ belongs to the class $\mathcal{L}_{\Delta_d}$ and
$$F^{*2}(x+\Delta_d)\sim 2F(x+\Delta_d),$$
then we say that the distribution $F$ belongs to the local subexponential distribution class $\mathcal{S}_{\Delta_d}$, where $F^{*n}$ is the $n$-th convolution of $F$ with itself for $n\in\mathbb{N}$ and $F^{*1}=F$. See Asmussen et al. \cite{AFK2003}. In particular, the classes $\mathcal{L}_{\Delta_\infty}$ and $\mathcal{S}_{\Delta_\infty}$ are called the long-tailed distribution class and subexponential distribution class, denoted by $\mathcal{L}$ and $\mathcal{S}$, respectively. At this time, we do not require $F\in\mathcal{L}$ in the definition of the subexponential distribution.

In the aforementioned two concepts, we replace ``for some $d$" with ``for each $d$", then we say that the distribution $F$ belongs to the local distribution class $\mathcal{L}_{loc}$ or $\mathcal{S}_{loc}$. See Borovkov and Borovkov \cite{BB2008}. Clearly, for each $d\in\mathbb{R}^+\cup\{\infty\}$, the following two inclusion relations that $\mathcal{L}_{loc}\subset\mathcal{L}_{\Delta_d}$ and $\mathcal{S}_{loc}\subset\mathcal{S}_{\Delta_d}$ are proper.

On the research concerning the local subexponential distribution, besides the above-mentioned papers, the readers can refer to Asmussen et al. \cite{AKKKT2002}, Ng and Tang \cite{NT2004}, Wang et al. \cite{WCW2005}, Shneer \cite{S2006}, Wang et al. \cite{WYWC2007}, Denisov and Shneer \cite{DS2007}, Denisov et al. \cite{DDS2008}, Cui et al. \cite{CWW2009}, Chen et al. \cite{CWC2009}, Yu et al. \cite{YWY2010}, Watanabe and Yamamuro \cite{WY2010}, Watanabe and Yamamuro \cite{WY2010-2}, Yang et al. \cite{YLS2010}, Wang and Wang \cite{WW2011}, Lin \cite{L2012}, Foss et al. \cite{FKZ2013}, Chen et al. \cite{CYW2013}, Cui et al. \cite{CWW2016}, Wang et al. \cite{WXCY2018}, Cui et al. \cite{CWX2018}, etc.

Among them, some results, for example, Lemma 2.3 and Corollary 2.1 of Denisov et al. \cite{DDS2008}, Proposition 6.1 of Wang and Wang \cite{WW2011} and Lemma 4.27 of Foss et al. \cite{FKZ2013}, required the following condition: let $F$ be a distribution in $\mathcal{S}_{\Delta_d}$ for some constant $d\in\mathbb{R}^+$, there exists a constant $x_0=x_0(F,d)\in\mathbb{R}^+\cup\{0\}$ such that $F(x+\Delta_{d})>0$ for all $x\ge x_0$ and

\begin{eqnarray}\label{101}
\sup_{x_0\le x\le y<\infty}\frac{F(y+\Delta_{d})}{F(x+\Delta_{d})}=:C=C(F,x_0)<\infty.
\end{eqnarray}
For the class $\mathcal{S}_{loc}$, it is required that, there is a constant $x_0=x_0(F)$ such that (\ref{101}) holds for all $d\in\mathbb{R}^+$. The condition is called that the local distribution of the distribution $F$ is almost decreasing, or the distribution is locally almost decreasing.
We know that, for many common distributions in the class $\mathcal{S}_{\Delta_d}$, their local distributions are almost decreasing with $C<\infty$, even $C=1$. Therefore, there also is a problem as follows.

\textbf{Problem 1.2.} Is each distribution in the class $\mathcal{S}_{\Delta_d}$ for some constant $d$, even in the class $\mathcal{S}_{loc}$, locally almost decreasing?

Based on Proposition 1.1 and Proposition 1.2, we have the following answer to the problem.
\begin{coro}\label{coro101}
For the class $\mathcal{S}_{loc}$ and the class $\mathcal{L}_{loc}\setminus\mathcal{S}_{loc}$ , there respectively exists a distribution $F$ such that its local distribution is not almost decreasing.
\end{coro}

We prove the above all results in the  next section. And in Section 3, we point out that, for the local subexponentiality of a distribution supported on $\mathbb{R}$, the local almost decrease is necessary in some sense. Thus, we find a substantial difference between the subexponential distribution and the local subexponential distribution.

%%%%%%%%%%%%%%%%%%%%%%%%%%%% section 7%%%%%%%%%%%%%%%%%%%%%%%%%%%
\setcounter{equation}{0}
\setcounter{lemma}{0}
\section{Proofs of the results}

{\bf Proof of Proposition \ref{pron101}}.
Let $\{a_n=2^{{n^2}},b_n=a_n+a_{m_n}\ln^2(n+1):n\in\mathbb{N}\}$ be a sequence of positive numbers, where $m_n=\min\{k:k\ge\sqrt{5}\big(\sqrt{6}\big)^{-1}n\}$. Clearly,
\begin{eqnarray}\label{500}
\sqrt{5}\big(\sqrt{6}\big)^{-1}n\le m_n<\sqrt{5}\big(\sqrt{6}\big)^{-1}(n+1)+1.
\end{eqnarray}
And let $f_0(\cdot):\mathbb{R}^+\cup\{0\}\longmapsto \mathbb{R}^+$ be a linear function defined as
$$f_0(0)=1,\ f_0(a_n)=f_0\big(\frac{a_{n+1}+b_n}{2}\big)=2a_n^{-3},\ f_0(b_n)=\frac{f_0(a_n)}{\ln(n+1)}$$
and
\begin{eqnarray}\label{501}
&&f_0(x)=\Big(f_0(0)+\frac{f_0(a_1)-f_0(0)}{a_1-a_0}x\Big)\textbf{1}(x\in[0, a_1])\nonumber\\
&&\ \ \ \ +
\sum_{n=1}^\infty\Big(\big(f_0(a_n)+\frac{f_0(b_n)-f_0(a_n)}{b_n-a_n}(x-a_n)\big)\textbf{1}(x\in J_{n1})\nonumber\\
&&\ \ \ \ +\big(f_0(b_n)+\frac{f_0(\frac{a_{n+1}+b_n}{2})-f_0(b_n)}{\frac{a_{n+1}+b_n}{2}-b_n}(x-b_n)\big)\textbf{1}(x\in J_{n2})\nonumber\\
&&\ \ \ \ +\big(f_0(\frac{a_{n+1}+b_n}{2})+\frac{f_0(a_{n+1})-f_0(\frac{a_{n+1}+b_n}{2})}
{\frac{a_{n+1}+b_n}{2}-b_n}(x-\frac{a_{n+1}+b_n}{2})\big)\textbf{1}(x\in J_{n3})\Big),
\end{eqnarray}
where $J_{n1}=(a_n,b_n],\ J_{n2}=(b_n,\frac{a_{n+1}+b_n}{2}]\ \text{and}\ J_{n3}=(\frac{a_{n+1}+b_n}{2},a_{n+1}]\ \text{for}\ n\in\mathbb{N}.$

Because
$\int_{a_n}^{a_{n+1}} f_0(y)dy\le f_0(a_n)a_{n+1}=2^{-2n^2+2(n+1)},$
$$0<a=\int_{0}^{\infty} f_0(x)dx<\infty.$$
Therefore, the function $f=a^{-1}f_0$ is a density corresponding to a distribution $F$ supported on $\mathbb{R}^+\cup\{0\}$. For the sake of simplicity, we set $a=1$ without affecting the results.

The density $f$ is not almost decreasing. In fact, when $n\to\infty$,
$$f(b_n)=f(a_n)\big(\ln(n+1)^{-1}\big)=o\big(f(a_n)\big)=o\big(f(\frac{a_{n+1}+b_n}{2})\big).$$

Then by (\ref{501}), when $n\to\infty$, we have
$$\frac{f(a_n)-f(b_n)}{b_n-a_n}=\frac{f(a_n)\big(1-\ln^{-1}(n+1)\big)}{a_{m_n}\ln^2(n+1)}=o\Big(\frac{f(b_n)}{a_{m_n}}\Big),$$
$$\frac{f(\frac{a_{n+1}+b_n}{2})-f(b_n)}{\frac{a_{n+1}+b_n}{2}-b_n}\sim\frac{2f(a_n)\big(1-\ln^{-1}(n+1)\big)}{a_{n+1}}
=o\Big(\frac{f(b_n)\ln^2(n+1)}{a_{n+1}}\Big),$$
and
$$\frac{f(\frac{a_{n+1}+b_n}{2})-f(a_{n+1})}{a_{n+1}-\frac{a_{n+1}+b_n}{2}}\sim\frac{2f(a_n)-f(a_{n+1})}{a_{n+1}}=o\big(f(a_{n+1})a^{-1}_nn^{5n}\big).$$
Thus,
\begin{eqnarray}\label{502}
\sup_{y\in J_{ni}}\mid f^{\prime}(y)\mid=o\big(\inf_{x\in J_{ni}}f(x)\big),\ i=1,2,3.
\end{eqnarray}

For $n\in\mathbb{N}^+$, we denote any two adjacent numbers in set $\{a_n,b_n,\frac{a_{n+1}+b_n}{2},a_{n+1}\}$ by $c_n$ and $d_n$. By the method of Lemma 4.1 in Xu et al. \cite{XFW2015} and (\ref{502}), we know that, for any fixed constant $t\in\mathbb{R}^+$ and variable $x\in\mathbb{R}^+$, there is a positive integer $n$ such that
\begin{eqnarray*}
f(x)&=&f(x-t)+\int_{x-t}^xf^{\prime}(y)dy\textbf{1}(c_n<x-t<x\le d_n)\nonumber\\
&&\ \ \ \ \ \ \ \ \ \ \ \ \ \ +\Big(\int_{x-t}^{d_n}+\int_{d_n}^x\Big)f^{\prime}(y)dy\textbf{1}(c_n<x-t\le d_n<x)\nonumber\\
&\sim& f(x-t).
\end{eqnarray*}
Thus, the density $f$ belongs to the class $\mathcal{L}_0$.

%In the following, we prove that
%$F\in\mathcal{S}_{loc}$. To this end, we only to prove that
%\begin{eqnarray}\label{502}
%\int_0^x f(x-y)f(y)dy\sim 2f(x).
%\end{eqnarray}
%\begin{eqnarray}\label{503}
%&&\int_0^{a_{n+1}} f(a_{n+1}-y)f(y)dy=\Big(\int_0^{a_n}+\int_{a_n}^{\frac{a_{n+1}+b_n}{2}}+\int_{\frac{a_{n+1}+b_n}{2}}^{a_{n+1}}\Big)
%f(a_{n+1}-y)f(y)dy\nonumber\\
%&&\ \ \ \ \ \ \ \ \ \ \ \ \ \ =\sum_{i=1}^3 I_i(a_{n+1})\sim 2f(a_{n+1}),\ as\ n\to\infty.
%\end{eqnarray}
%From (\ref{501}), we know that, when $n\to\infty$,
%$$f(a_{n+1}-a_n)=f\big(\frac{a_{n+1}+b_n}{2}\big)+\frac{f(a_{n+1})-f(\frac{a_{n+1}+b_n}{2})}{\frac{a_{n+1}-b_n}{2}}\cdot\frac{a_{n+1}-a_n-b_n}{2}
%\sim f(a_{n+1}).$$
%And then by the following two facts that
%\begin{eqnarray*}
%I_1(a_{n+1})\le f(a_{n+1}-a_n)\int_0^{a_n}f(y)dy\sim f(a_{n+1}),\ as\ n\to\infty
%\end{eqnarray*}
%and
%\begin{eqnarray*}
%I_1(a_{n+1})\ge f(a_{n+1})\int_0^{a_n}f(y)dy\sim f(a_{n+1}),\ as\ n\to\infty,
%\end{eqnarray*}
%we have
%\begin{eqnarray}\label{504}
%I_1(a_{n+1})\sim f(a_{n+1}).
%\end{eqnarray}
%Since $\frac{a_{n+1}-b_n}{2}\le a_{n+1}-y<a_{n+1}-a_n$ for $a_n<y\le \frac{a_{n+1}+b_n}{2}$,
%\begin{eqnarray}\label{505}
%I_2(a_{n+1})\le \big(f(a_n)\big)^2a_{n+1}=o\big(f(a_{n+1})\big),\ as\ n\to\infty.
%\end{eqnarray}
%By (\ref{504}) and (\ref{505}), we have
%\begin{eqnarray}\label{506}
%&&I_3(a_{n+1})=\Big(\int_0^{a_n}+\int_{a_n}^{\frac{a_{n+1}+b_n}{2}}\Big)f(a_{n+1}-y)f(y)dy\sim f(a_{n+1}),\ as\ n\to\infty.
%\end{eqnarray}
%Therefore, according to (\ref{504}), (\ref{505}) and (\ref{506}), (\ref{503}) is holds.

Now, we go on to prove that $f\in\mathcal{S}_0$, that is for $x\in J_n=(a_n,a_{n+1}]$,
\begin{eqnarray}\label{507}
I(x)&=&\int_0^x f(x-y)f(y)dy=\Big(2\int_0^{a_{m_n}}+ \int_{a_{m_n}}^{x-a_{m_n}}\Big) f(x-y)f(y)dy\nonumber\\
&=&2I_1(x)+I_2(x)\sim 2f(x).
\end{eqnarray}

For $ x\in J_n$, because $y$ and $x-y\in(a_{m_n},x-a_{m_n}]$, by (\ref{501}), we have
\begin{eqnarray}\label{508}
I_2(x)&\le& f^2(a_{m_n})a_{n+1}=2^{-6m_n^2+(n+1)^2+2}
\le 2^{-3(n+1)^2+3(n+1)^2-5n^2+(n+1)^2+2}\nonumber\\
&=&o\big(f(a_{n+1})\big)=o\big(f(x)\big).
\end{eqnarray}

In the following, we deal with $I_1(x)$. Because $f\in\mathcal{L}_0$, we just have to prove that
\begin{eqnarray}\label{5080}
f(x-a_{m_n})\sim f(x)
\end{eqnarray}
for $x\in J_{ni},i=1,2,3$, respectively.

When $x\in J_{n1}=(a_n,b_n]$, because $0\le y\le a_{m_n}$,
$$\frac{a_{n}+b_{n-1}}{2}<a_n-a_{m_n}<x-a_{m_n}\le x-y\le x\le b_n.$$
If $a_n\le x-a_{m_n}<x\le b_n$, then by (\ref{501}), $f(b_n)=\frac{f(a_n)}{\ln(n+1)}$ and (\ref{500}),
\begin{eqnarray*}
f(x-a_{m_n})&=&f(x)+\frac{f(a_n)\big(1-\ln^{-1}(n+1)\big)}{a_{m_n}\ln^2(n+1)}\nonumber\\
&=&f(x)+\frac{f(b_n)\big(1-\ln^{-1}(n+1)\big)}{a_{m_n}\ln(n+1)}\sim f(x).
\end{eqnarray*}
If $\frac{a_{n}+b_{n-1}}{2}< a_n-a_{m_n}\le x-a_{m_n}<a_n<x\le b_n$, then by (\ref{501}) and (\ref{500}),
\begin{eqnarray*}
f(a_n)&\le& f(x-a_{m_n})\le f(a_n-a_{m_n})\nonumber\\
&=&f(a_{n})+\frac{2\big(f(a_{n-1})-f(a_n)\big)a_{m_n}}{a_n-b_{n-1}}\nonumber\\
&\le&f(a_{n})+f(a_{n})2^{-\frac{1}{6}n^2+6n}\sim f(a_n),
\end{eqnarray*}
that is $f(x-a_{m_n})\sim f(a_n)$. On the other hand, because $x-a_n\le a_{m_n}$, we have
\begin{eqnarray*}
f(x)=f(a_n)+\frac{f(a_n)\big(\ln^{-1}(n+1)-1\big)}{a_{m_n}\ln^2(n+1)}\cdot (x-a_n)\sim f(a_n),
\end{eqnarray*}
that is (\ref{5080}) holds for $x\in J_{n1}=(a_n,b_n]$.

When $x\in J_{n2}=(b_n,\frac{a_{n+1}+b_{n}}{2}]$, then
$$a_{n}<b_n-a_{m_n}<x-a_{m_n}<x\le \frac{a_{n+1}+b_{n}}{2}.$$
If $b_{n}<x-a_{m_n}<x\le \frac{a_{n+1}+b_{n}}{2},$ then by (\ref{501}) and (\ref{500}), we have
$$f(x-a_{m_n})=f(x)-\frac{2\big(f(a_n)-f(b_n)\big)}{a_{n+1}-b_n}\cdot a_{m_n}\sim f(x).$$
If $a_n<b_n-a_{m_n}<x-a_{m_n}\le b_n<x,$ then by (\ref{501}), we have
$$f(x-a_{m_n})\le f(b_n-a_{m_n})=f(b_n)+\frac{f(a_n)\big(1-\ln^{-1}(n+1)\big)}{a_{m_n}\ln^2(n+1)}\cdot a_{m_n}\sim f(b_n).$$
On the other hand, because $x-b_n\le a_{m_n}$, we have
\begin{eqnarray*}
f(x)=f(b_n)+\frac{2\big(f(a_n)-f(b_n)\big)}{a_{n+1}-b_n}\cdot (x-b_n)\sim f(b_n)\le f(x-a_{m_n}).
\end{eqnarray*}
Thus, (\ref{5080}) holds for $x\in(b_n,\frac{a_{n+1}+b_{n}}{2}]$.

When $x\in J_{n3}=(\frac{a_{n+1}+b_{n}}{2},a_{n+1}]$, then
$$\frac{a_{n+1}+b_{n}}{2}-a_{m_n}<x-a_{m_n}<x\le a_{n+1}.$$
If $\frac{a_{n+1}+b_{n}}{2}<x-a_{m_n}<x\le a_{n+1}$, then
$$f(x-a_{m_n})=f(x)+\frac{2\big(f(a_n)-f(a_{n+1})\big)}{a_{n+1}-b_n}\cdot a_{m_n}\sim f(x).$$
If $b_n<\frac{a_{n+1}+b_{n}}{2}-a_{m_n}<x-a_{m_n}\le\frac{a_{n+1}+b_{n}}{2}<x\le a_{n+1}$, then by (\ref{501}),
$$f(x-a_{m_n})\ge f(\frac{a_{n+1}+b_{n}}{2}-a_{m_n})=f(a_n)+\frac{2\big(f(a_n)-f(b_n)\big)}{a_{n+1}-b_n}\cdot a_{m_n}\sim f(a_n);$$
and by $x-\frac{a_{n+1}+b_{n}}{2}\le a_{m_n}$,
$$f(x)=f(a_n)+\frac{2\big(f(a_n)-f(a_{n+1})\big)}{a_{n+1}-b_n}\cdot \big(x-\frac{a_{n+1}+b_{n}}{2}\big)\sim f(a_n)\ge f(x-a_{m_n}).$$
Thus, (\ref{5080}) holds for $x\in(\frac{a_{n+1}+b_{n}}{2},a_{n+1}]$.

Therefore, (\ref{507}) is proved , that is $f\in\mathcal{S}_0$.\hfill$\Box$
\\

{\bf Proof of Proposition \ref{pron102}.} Cline \cite{C1986} gave a distribution $F_0$ on $\mathbb{R}^+\cup\{0\}$ such that
$$\overline{F_0}(x)=e^{-\alpha x-\chi(x)}\textbf{1}\big(x\in\mathbb{R}^+\cup\{0\}\big)+\textbf{1}\big(x\in\mathbb{R}^-\big)\ \ \text{for}\ x\in\mathbb{R},$$
where $\chi(x)=x^{2^{-1}+\delta\cos\circ\ln (x+1)}\textbf{1}\big(x\in\mathbb{R}^+\cup\{0\}\big)$, $\alpha\in\mathbb{R}^+$ and $\delta\in(0,2^{-1})$.
Then $F_0\in\mathcal{L}(\alpha)\setminus\mathcal{S}(\alpha)$, where the two distribution classes are defined as follows:
$$\mathcal{L}(\alpha)=\{V\ \text{on}\ \mathbb{R}^+\cup\{0\}: \overline{V}(x+t)\sim e^{-\alpha t}\overline{V}(x)\ \text{for each constant}\ t\in\mathbb{R}\}$$
and
$$\mathcal{S}(\alpha)=\Big\{V\ \text{in}\ \mathcal{L}(\alpha): M(V,\alpha)=\int_0^\infty e^{\alpha y}V(dy)<\infty\ \text{and}\ \overline{V^{*2}}(x)\sim2M(V,\alpha)\overline{V}(x)\Big\}.$$

Further, we define a density
$$f(x)=ae^{\alpha x}\overline{F_0}(x)\textbf{1}(x\in\mathbb{R}^+\cup\{0\})=ae^{-\chi(x)}\textbf{1}(x\in\mathbb{R}^+\cup\{0\})\ \ \text{for}\ x\in\mathbb{R},$$
where $0<a^{-1}=\int_0^\infty e^{-\chi(y)}dy<\infty$. Clearly, $f\in\mathcal{L}_0$. Then we define a distribution $F_0^I$ such that
$$\overline{F_0^I}(x)=b\int_x^\infty\overline{F_0}(y)dy\textbf{1}(x\in\mathbb{R}^+\cup\{0\})
+\textbf{1}\big(x\in\mathbb{R}^-\big)$$
for $x\in\mathbb{R}$, where $0<b^{-1}=\int_0^\infty \overline{F_0}(y)dy<\infty$. By Karamata's theorem, $\overline{F_0^I}(x)\sim \alpha^{-1}b\overline{F_0}(x)$. Thus, $F_0^I\in\mathcal{L}(\alpha)\setminus\mathcal{S}(\alpha)$ and there is a positive function $h(\cdot)$ such that $h(x)\uparrow\infty,\ x^{-1}h(x)\to0$ and $\overline{F_0^I}(x-y)\sim e^{\alpha y}\overline{F_0^I}(x)$ uniformly for all $\mid y \mid\le h(x)$. For the function, we have
\begin{eqnarray}\label{201}
&&\int_{h(x)}^{x-h(x)}f(x-y)f(y)dy\big(f(x)\big)^{-1}=a\int_{h(x)}^{x-h(x)}\overline{F_0}(x-y)\overline{F_0}(y)dy\big(\overline{F_0}(x)\big)^{-1}\nonumber\\
&\sim&a\alpha b^{-1}\int_{h(x)}^{x-h(x)}\overline{F_0^I}(x-y)\overline{F_0^I}(dy)\big(\overline{F_0^I}(x)\big)^{-1}.
\end{eqnarray}
By (\ref{201}), $F_0^I\notin\mathcal{S}(\alpha)$ and Proposition 6 of Asmussen et al. \cite{AFK2003}, $f\notin\mathcal{S}_0$.

Finally, we prove that $f$ is not almost decreasing. In fact, the conclusion follows from
$$f(e^{2k\pi-1})\big(f(e^{(2k+1)\pi-1})\big)^{-1}\to\infty,\ \ \text{as}\ k\to\infty.$$

In this way, we have completed the proof of Proposition \ref{pron102}.\hfill$\Box$\\

{\bf Proof of Corollary \ref{coro101}.} The two conclusions directly follow from the fact
$$F(x+\Delta_d)\sim df(x).$$ \hfill$\Box$
\\

%%%%%%%%%%%%%%%%%%%%%%%%%%%% section 8%%%%%%%%%%%%%%%%%%%%%%%%%%%
\setcounter{equation}{0}
\setcounter{lemma}{0}
\section{The role of the almost decrease}

In this section, we illustrate the role of the local almost decrease of a distribution.

For a distribution $F$ on $\mathbb{R}$, similar to the definition of subexponential distribution, there are two different definitions of the local subexponential distribution.

A distribution $F$ on $\mathbb{R}$ belongs to the local subexponential distribution class $\mathcal{S}_{\Delta_d}$ for some constant $d\in\mathbb{R}^+\cup\{\infty\}$, if $F\in\mathcal{L}_{\Delta_d}$ and
\begin{eqnarray}\label{301}
F^{*2}(x+\Delta_d)\sim 2F(x+\Delta_d).
\end{eqnarray}
See (4.42) of Foss et al. \cite{FKZ2013} or Definition 1.3 (iii) of Watanabe and Yamamuro \cite{WY2010}.

Before giving the second definition, we introduce some concepts and notations. Let $X$ be a r.v. with distribution $F$ on $\mathbb{R}$ which means $q_1:=F[0,\infty)>0$ and $q_2:=F(-\infty,0)<1$. Two distributions corresponding to r.v.s $X_1$ and $-X_2$ are respectively defined by $F_1$ on $\mathbb{R^+}\cup\{0\}$ and $F_2$ on $\mathbb{R^-}$ such that
$$F_1(dy)=q_1^{-1}F(dy)\textbf{1}_{\mathbb{R^+}\cup\{0\}}(y)\ \text{and}\ F_2(dy)=q_2^{-1}F(dy)\textbf{1}_{\mathbb{R^-}}(y)$$
for all $y\in\mathbb{R}$. Then
\begin{eqnarray}\label{302}
F=q_1F_1+q_2F_2
\end{eqnarray}
and for all $x\in\mathbb{R}$ and some $0<d<\infty$,
\begin{eqnarray}\label{303}
F^{*2}(x+\triangle_d)&=&q_1^2F_1^{*2}(x+\triangle_d)+2q_1q_2F_1*F_2(x+\triangle_d).
\end{eqnarray}

The second definition is as follows. A distribution $F$ supported on $\mathbb{R}$ belongs to the local subexponential distribution class $\mathcal{S}_{\Delta_d}$ for some constant $0<d\le\infty$, if $F_1\in\mathcal{L}_{\Delta_d}$ and
\begin{eqnarray}\label{304}
F_1^{*2}(x+\Delta_d)\sim 2F_1(x+\Delta_d).
\end{eqnarray}
See Section 4.2 in Foss et al. \cite{FKZ2013} or Subsection 1.1 in Wang and Wang \cite{WW2011}.

Similarly, we can also give another two definitions for $F$ on $\mathbb{R}$ belonging to $\mathcal{S}_{loc}$.
Therefore, we are naturally concerned with the following problem.

\textbf{Problem 3.1.}  When $F$ on $\mathbb{R}$ belongs to the class $\mathcal{L}_{\Delta_d}$, can the two representations (\ref{301}) and (\ref{304}) be implied from each other?

%Clearly, if $F_1\in {\cal S}_{\triangle_d}$, then $F\in {\cal L}_{\triangle_d}$.
Lemma 4.26 of Foss et al. \cite{FKZ2013} show that (\ref{304}) can follow from (\ref{301}). Conversely,
%that is $F_1\in{\cal L}_\triangle$ and \begin{eqnarray}\label{504} F_1^{*2}(x+\triangle)&\sim&2F_1(x+\triangle)=2 q_1^{-1}F(x+\triangle), \end{eqnarray}
under one of the following two conditions that $F_1$ is locally almost decreasing and $X_2$ is a light-tailed r.v. that there is a constant $\varepsilon\in\mathbb{R}^+$ such that $Ee^{\varepsilon X_2}<\infty$, then (\ref{301}) can follow from (\ref{304}). See Lemma 4.27 of Foss et al. \cite{FKZ2013} for the previous condition and Proposition 6.1 of Wang and Wang \cite{WW2011} for both.
%\begin{eqnarray}\label{505}
%F_1*F_2(x+\triangle_d)=\int_{-\infty}^0F_1(x-y+\triangle_d)dF_2(y)\sim F_1(x+\triangle_d)=q_1^{-1}F(x+\triangle_d).
%\end{eqnarray}
%Thus, by (\ref{502}), we have
%\begin{eqnarray}\label{503}
%F^{*2}(x+\triangle_d)&\sim&2F(x+\triangle_d).
%\end{eqnarray}

As Wang and Wang \cite{WW2011} point out that, however, ``we do not know whether or not these conditions can be cancelled". We also have not found any relevant discussions on this issue in other references. There we give a answers to the question, thus to Problem 3.1.

\begin{pron}\label{pron301}
There exists a distribution $F$ on $\mathbb{R}$ such that, $F_2$ is heavy-tailed, $F_1$ is not almost decreasing and belongs to the class $\mathcal{S}_{loc}$, that is (\ref{304}) holds, but (\ref{301}) does not hold, that is $F$ not belong to the same one.
\end{pron}
In other words, if (\ref{304}) holds, then in order to get (\ref{301}), the two conditions that, the distribution $F_1$ is locally almost decreasing and the distribution $F_2$ is light-tailed, are necessary in a sense.\\
\\
{\bf Proof of Proposition \ref{pron301}.}
We assume that the distribution $F_1$ in (\ref{302}) has a density $f$ as (\ref{501}), then $F_1\in {\cal S}_{loc}$ and is not locally almost decreasing. Then we have
\begin{eqnarray}\label{305}
F_1*F_2(x+\triangle)&=&\int_{-\infty}^0F_1(x-y+\triangle_d)dF_2(y)\nonumber\\
&\sim&d\int_{-\infty}^0f(x-y)dF_2(y)%=q_1^{-1}d\int_{\infty}^xf(z)dP(-X_2\le x-z)
\nonumber\\
&=&d\int_x^{\infty}f(z)dP(X_2\le z-x).
\end{eqnarray}
In (\ref{305}), we take $x=b_n$ and denote $c_n=\frac{a_{n+1}+b_n}{2}$ for each $n\in\mathbb{N}$, then
\begin{eqnarray}\label{306}
F_1*F_2(b_n+\triangle)&\sim&d\sum_{m=n}^\infty\int_{b_m}^{b_{m+1}}f(z)dP(X_2\le z-b_n)\nonumber\\
&\sim&d\int_{b_n}^{b_{n+1}}f(z)dP(X_2\le z-b_n),
\end{eqnarray}
the last step in (\ref{306}) comes from the following two facts that
\begin{eqnarray*}
\sum_{m=n+1}^\infty\int_{b_m}^{b_{m+1}}f(z)dP(X_2\le z-b_n)&\le& \sum_{m=n+1}^\infty f(b_n)\frac{f(a_{m})}{f(b_n)}\int_{b_m}^{b_{m+1}}dP(X_2\le z-b_n)\nonumber\\
&\le&f(b_n)\frac{f(a_{n+1})}{f(b_n)}P(X_2> b_{n+1}-b_n)
\end{eqnarray*}
and
\begin{eqnarray*}
I(n)&:=&\int_{b_n}^{b_{n+1}}f(z)dP(X_2\le z-b_n)\ge\int_{b_n}^{c_{n}}f(z)dP(X_2\le z-b_n)\nonumber\\
&\ge& f(b_n)P(0<X_2\le c_{n}-b_n)\sim f(b_n).
\end{eqnarray*}

We further analyze $I(n)$. Because $f$ is a linear function defined by (\ref{501}), for each $n\in\mathbb{N}$, we can take $b_n\le d_n\le c_n$ and $c_n\le s_n\le a_{n+1}$ such that
$$f(d_n)=f(s_n)=f(b_n)\sqrt{\ln(n+1)}=\frac{f(a_n)}{\sqrt{\ln(n+1)}}=\frac{2}{a_n^3\sqrt{\ln(n+1)}},$$
then $d_n\le y+b_n\le s_n$ for $d_n-b_n\le y\le s_n-b_n$ and
\begin{eqnarray*}
I(n)&\ge&\int_{d_n-b_n}^{s_n-b_n}f(y+b_n)dP(X_2\le y)\nonumber\\
&\ge&f(b_n)\sqrt{\ln(n+1)}P(d_n-b_n<X_2\le s_n-b_n).
\end{eqnarray*}
Now, we define a concrete distribution $F_2$ of r.v. $-X_2$ such that
\begin{eqnarray*}
F_2(b_{n_m}-s_{n_m},b_{n_m}-d_{n_m}]=P(d_{n_m}-b_{n_m}<X_2\le s_{n_m}-b_{n_m})=c\big(\ln(n_m+1)\big)^{-\frac{1}{3}}
\end{eqnarray*}
for a subsequence of positive integers $\{n_m:\ m\in\mathbb{N}\}$, for example, we can take $n_m=[e^{m^4}]$, where $c>0$ is a regularized constant. Clearly,
\begin{eqnarray*}
\frac{I(n_m)}{f(b_{n_m})}\to \infty\ \ \text{as}\ m\to\infty.
\end{eqnarray*}

Therefore, for each $d\in\mathbb{R}^+$,
\begin{eqnarray*}
\frac{F_1*F_2(b_{n_m}+\triangle_d)}{F(b_{n_m}+\triangle_d)}\sim\frac{F_1*F_2(b_{n_m}+\triangle_d)}{f(b_{n_m})}\to \infty\ \ \text{as}\ m\to\infty,
\end{eqnarray*}
which means
\begin{eqnarray}\label{307}
\limsup\frac{F^{*2}(x+\triangle_d)}{F(x+\triangle_d)}=\infty.
\end{eqnarray}
In addition, the distribution $F_2$ is heavy-tailed, otherwise, by Proposition 6.1 of Wang and Wang \cite{WW2011}, (\ref{301}) holds, which is conflicting with (\ref{307}). $\hfill\Box$
\vspace{0.2cm}

We already know that, when $d=\infty$, $F_1\in\mathcal{S}$ is equivalent to $F\in\mathcal{S}$. However, when $d<\infty$, $F_1\in\mathcal{S}_{loc}$ can not guarantee $F\in\mathcal{S}_{loc}$. This is a substantial difference between the subexponential distribution and the local subexponential distribution.

\vspace{0.8cm}

\textbf{Acknowledgements} The authors are very grateful to Dr. Changjun Yu and doctoral student Hui Xu for their helpful discussions and comments.


\begin{thebibliography}{99}

%\bibitem{Albrecher2006} Albrecher, H., Teugels, J. L., 2006. Exponential behavior in the presence of dependence in risk theory. Journal of  Applied Probability 43, 257-273.

%\bibitem {Asimit2010} Asimit, A., Badescu, A. L., 2010. Extremes on the discounted aggregate claims in a time dependent risk model. Scandinavian Actuarial Journal 2, 93-104.
\bibitem{AKKKT2002} Asmussen, S., Kalashnikov, V., Konstantinides, D., Kl\"{u}ppelberg, C.
and Tsitsiashvili, G., 2002. A local limit theorem for random walk
maxima with heavy tails. Statist. Probab. Lett. 56, 399-404.


\bibitem{AFK2003} Asmussen, S., Foss, S., Korshunov, D., 2003. Asymptotics for sums of
random variables with local subexponential behavior. J. Theor.
Probab. 16, 489-518.

%\bibitem{Asmussen2002_K} Asmussen, S., Kalashnikov, V., Konstantinides, D., Kl\"{u}ppelberg, C., Tsitsiashvili, G., 2002. A local limit theorem for random walk maxima with heavy tails. Statistics \& Probability Letters 56, 399-404.

%\bibitem{Badescu2009} Badescu, A. L., Cheung, E. C. K., Landriault, D., 2009. Dependent risk models with bivariate phase-type distributions. Journal of Applied Probability 46, 113-131.

\bibitem{BB2008} Borovkov, A.A., Borovkov, K.A., 2008. Asymptotic Analysis of random walks. Cambridge: Cambridge University Press.

\bibitem{BGT1987} Bingham, N. H.,  Goldie, C. M., Teugels, J. L., 1987. Regular Variation. Cambridge: Cambridge University Press.

%\bibitem{CYZ2003} Chan, W., Yang, H., Zhang, L., 2003. Some results on ruin probabilities in a two-dimensional risk model. Insurance: Mathematics \& Economics 32, 345-358.

%\bibitem{} Chen, Y., Chen, A.,  Ng, K. W., 2010. The strong law of large numbers for extended negatively dependent random variables. Journal of  Applied Probability 47, 908-922.

%\bibitem{CN2007} Chen, Y., Ng, K. W., 2007. The ruin probability of the renewal model with constant interest force and negatively dependent heavy-tail claims. Insurance: Mathematics \& Economics 40, 415-423.

\bibitem{CYW2013} Chen, W., Yu, C., Wang, Y., 2013. Some discussions on the local distribution classes. Statistics \& Probability Letters 83, 1654-1661.

%\bibitem{Chen2013} Chen, Y., Wang, Y., Wang, K., 2013. Asymptotic results for ruin probability of a two-dimensional renewal risk model. Stochastic Analysis and Applications 31, 80-91.


%\bibitem{CYN2011} Chen, Y., Yuen, K. C., Ng, K. W., 2011. Asymptotics for the ruin probabilities of a two-dimensional renewal risk model with heavy-tailed claim. Applied Stochastic Models in Business and Industry 27, 290-300.

\bibitem{CWC2009} Chen, G., Wang, Y., Cheng, F., 2009. The uniform local asymptotics
of the overshoot of a random walk with heavy-tailed increments. Stochastic Models 25, 508-521.

\bibitem{CNW1973} Chover, J., Ney, P., Wainger, S., 1973. Functions of probability measures. J. Anal. Math. 26, 255¨C302.

\bibitem{C1986} Cline, D. B. H., 1986. Convolution tails, product tails and Domains of attraction. Probab. Theor. Rel. Fields 72, 529-557.

\bibitem{CWX2018} Cui Z., Wang, Y., Xu, H., 2018. Some positive conclusions related to the Embrechts-Goldie conjecture. Submited to Stochastics: An International Journal Of Probability And Stochastic Processes.

\bibitem{CWW2009} Cui Z., Wang, Y., Wang, K., 2009. Asymptotics for the moments
of the overshoot and undershoot of a random walk. Advance in Applied Probability 41, 469-494.

\bibitem{CWW2016} Cui Z., Wang, Y., Wang, K., 2016. The uniform local asymptotics for a L¨¦vy process
and its overshoot and undershoot. Communications in Statistics-Theory and Methods, 45(4), 1156-1181.

\bibitem{DS2007} Denisov, D., Shneer, V., 2007. Local asymptotics of
the cycle maximum of a heavy-tailed random walk. Advance in Applied Probability 39, 221-244.

\bibitem{DDS2008} Denisov, D., Dieker, A. B., Shneer, V., 2008. Large deviations for random walks under subexponentiality: the big domain.
Ann. Probab. 36 (5), 1946-1991.

\bibitem{FKZ2013} Foss, S., Korshunov, D., Zachary, S., 2013.  An Introduction to
Heavy-tailed and Subexponential Distributions. Springer, Second Edition.

%\bibitem{FK2007} Foss, S., Korshunov, D., 2007. Lower limits and equivalences for convolution tails. Ann. Probab. 35, 366¨C383.
%\bibitem{Dong2011} Dong, Y., Wang, Y., 2011. Uniform estimates for ruin probabilities in the renewal risk model with upper-tail independent claims and premiums. Journal of Industrial and Management Optimization 7, 849-874.

%\bibitem{EmbKluMik1997} Embrechts, P., Kl\"{u}ppelberg, C. and Mikosch, T., 1997. Modelling Extremal Events for Insurance and Finance. Springer.

%\bibitem{} Gao, Q., Wang, Y., 2009. Ruin probability and local ruin probability in the random multi-delayed renewal risk model. Statistics \& Probability Letters 79(5), 588-596.

%\bibitem{HT2008} Hao, X., Tang, Q., 2008. A uniform asymptotic estimate for discounted aggregrate claims with subexponential tail. Insurance: Mathematics \& Economics 43, 116-120.


%\bibitem{HWZ2014} Huang, W,  Weng, C., Zhang, Y., 2014. Multivariate risk models under heavy-tailed risks. Applied Stochastic Models in Business and Industry 30, 341-360.


%\bibitem{KK2000} Kalashnikov, V.,  Konstantinides, D., 2000. Ruin under interest force and subexponential claims: a simple treatment. Insurance: Mathematics \& Economics 27, 145-149.

\bibitem{K1989} Kl\"{u}ppelberg, C., 1989. Subexponential distributions and characterizations of related classes. Probab. Theory Relat. Fields 82, 259-269.

\bibitem{K2006} Korshunov, D., 2006. On the distribution density of the supremum of a random walk in the subexponential case. Siberian Mathematical Journal, 47(6), 1060¨C1065.

%\bibitem{KTT2002} Konstantinides, D., Tang, Q., Tsitsiashvili, G., 2002. Estimates for the ruin probability in the classical risk model with constant interest force in the presence of heavy tail. Insurance: Mathematics \& Economics 31, 447-460.

%\bibitem{KNJ2000} Kotz, S., Balakrishnan, N. and Johnson, N. L., 2000. Continuous Multivariate Distributions. Volume 1: Models and Applications. Wiley.

%\bibitem{JWCX2015} Jiang, T., Wang, Y., Chen, Y., Xu, H., 2015. Uniform asymptotic estimate for finite-time ruin probabilities of a time-dependent bidimensional renewal model. Insurance: Mathematics \& Economics 64, 45-53.

%\bibitem{LLT2007} Li, J., Liu, Z., Tang, Q., 2007. On the ruin probabilities of a bidimensional perturbed risk model. Insurance: Mathematics \& Economics 41, 85-195.


%\bibitem{Li2010} Li, J., Tang, Q., Wu, R., 2010. Subexponential tails of discounted aggregate claims in a time-dependent renewal risk model. Advance in Applied Probability 42, 1126-1146.

\bibitem{L2012} Lin, J., 2012. Second order Subexponential Distributions with Finite Mean and Their Applications to Subordinated Distributions. J. Theor. Probab. 25, 834-853.

%\bibitem{LGW2012} Liu, X.,  Gao, Q., Wang, Y., 2012. A note on a dependent risk model with constant interest rate. Statistics \& Probability Letters 82, 707-712.

%\bibitem{N2006} Nelsen, R. B., 2006. An introduction to Copulas (Second edition). New York: Springer.

\bibitem{NT2004} Ng, K., Tang, Q., 2004. Asymptotic Behavior of Tail and Local
Probabilities for Sums of Subexponential Random Variables. Journal of Applied Probability 41, 108-116.

\bibitem{S1996} Sgibnev, S. M., 1996. On the distribution of the maxima of partial sums. Statistics \& Probability Letters 28, 235-238,

\bibitem{S2006} Shneer, V. V., 2006. Estimates for the interval probabilities of sums of random variables
with locally subexponential distributions. Siberian Math. J. 47, 946-955.

%\bibitem{} Tang, Q., Tsitsiashvili, G., 2003. Precise estimates for the ruin probability in finite horizon in a discrete-time model with heavy-tailed insurance and financial risks. Stochastic Processes and their Applications 108, 299-325.

%\bibitem{T2005} Tang, Q., 2005. The finite-time ruin probability of the compound Poisson model with constant interest force. Journal of Applied Probability 42, 608-619.

%\bibitem{T2007} Tang, Q., 2007. Heavy tails of discounted aggregrate claims in the continuous-time renewal model. Journal of  Applied Probability 44, 285-294.

%\bibitem{WWG2013} Wang, K., Wang, Y., Gao, Q., 2013. Uniform asymptotics for the finite-time ruin probability of a new dependent risk model with a constant interest rate. Methodology and Computing in Applied Probability 15, 109-124.

\bibitem{WCW2005} Wang, Y., Cheng, D., Wang, K., 2005. The closure of a local subexponential
distribution class under convolution roots with applications to the compound
Poisson process. Journal of Applied Probability 42, 1194-1203.

\bibitem{WW2006} Wang, Y., Wang, K., 2006. Asymptotics of the density of the supremum of a random walk with heavy-tailed increments.
Journal of Applied Probability 3, 874-879.

\bibitem{WW2011} Wang, Y., Wang, K., 2011. Random walks with non-convolution equivalent increments and their applications.
Journal of Mathematical Analysis and Applications 374, 88-105.

\bibitem{WXCY2018} Wang, Y., Xu, H., Cheng, D., Yu, C., 2018. The local asymptotic estimation for the supremum of a random walk, submitted for publication. Statistical Papers, 59, 99-126.

\bibitem{WYWC2007} Wang, Y., Yang, Y., Wang, K., Cheng, D., 2007.
Some new equivalent conditions on asymptotics and local asymptotics
for random sums and their applications. Insurance: Mathematics \& Economics
42(2), 256-266.

\bibitem{WY2010} Watanabe, T., Yamamura, K., 2010. Local Subexponentiality and Self-decomposability.
J. Theor Probab. 23(4),1039-1067

\bibitem{WY2010-2} Watanabe, T., Yamamura, K., 2010. Ratio of the tail of an infinitely divisible
distribution on the line to that of its L\'evy measure. Electron. J. Probab. 15(2), 44-74.

\bibitem{WY2017} Watanabe, T., Yamamura, K., 2017. Two non-closure properties on the class
of subexponential densities. J. Theor. Probab. 30(3),1059-1075.

\bibitem{XFW2015} Xu, H., Foss, S., Wang, Y., 2015.  On closedness under convolution and convolution roots of the class of long-tailed distributions. Extremes, 18, 605-628.

%\bibitem{YL2014} Yang, H., Li, J., 2014. Asymptotic finite-time ruin probability for a bidmensional renewal risk model with constant interest force and dependent subexponeantial claims. Insurance: Mathematics \& Economics 58, 185-192.

\bibitem{YLS2010} Yang, Y., Leipus, R., Siaulys, J., 2010. Local presice deviations for sums of random 
variables with O-reqularly varying densities. Statistics \& Probability Letters 80, 1559-1567.

%\bibitem{YW2013} Yang, Y., Wang, Y., 2013. Tail behavior of the product of two dependent random variables with applications to risk theory. Extremes 16, 55-74.

\bibitem{YWY2010} Yu, C., Wang, Y., Yang, Y., 2010. The closure of
the convolution equivalent distribution class under distribution
class with applications to random sums. Statistics \& Probability Letters 80, 462-472.


%\bibitem{} Yu, C., Wang, Y., Cheng, D., 2015. Tail behavior of the sums of dependent and  heavy-tailed  random variables. Journal of the Korean Statistical Society 44, 12-27.

%\bibitem{YGW2006} Yuen, K. C., Guo, J., Wu, X., 2006. On the first time of ruin in the bivariate compound Poisson model. Insurance: Mathematics \& Economics 38, 298-308.

\end{thebibliography}
\end{document}